\documentclass[oneside,english]{amsart}
\usepackage[T1]{fontenc}
\usepackage[latin9]{inputenc}
\usepackage{amsthm}
\usepackage{amssymb}
\PassOptionsToPackage{normalem}{ulem}
\usepackage{ulem}

\makeatletter

\providecommand{\tabularnewline}{\\}

\numberwithin{equation}{section}
\numberwithin{figure}{section}
\theoremstyle{plain}
\newtheorem{thm}{\protect\theoremname}
\theoremstyle{plain}
\newtheorem{lem}[thm]{\protect\lemmaname}

\makeatother

\usepackage{babel}
\providecommand{\lemmaname}{Lemma}
\providecommand{\theoremname}{Theorem}

\begin{document}
\title{Recurrent Relations for Multiple of Triangular Numbers being Triangular
Numbers}
\author{Vladimir PLETSER}
\address{European Space Agency (ret.) }
\address{Email: Pletservladimir@gmail.com}
\begin{abstract}
We search for triangular numbers that are multiples of other triangular
numbers. It is found that for any positive non-square integer multiplier,
there is an infinity of multiples of triangular numbers that are triangular
numbers and recurrent relations are deduced theoretically. If the
multiplier is a squared integer, there is either one or no solution,
depending on the multiplier value.
\end{abstract}

\maketitle

\section{Introduction}

Triangular numbers $T_{t}=\frac{t\left(t+1\right)}{2}$ enjoy many
interesting properties and many relations and formulas are associated
to them (see e.g. \cite{key-1,key-2}). We are interested in multiples
of triangular number $kT_{t}$ being triangular numbers
\begin{equation}
T_{\xi}=kT_{t}\label{eq:1}
\end{equation}
Euler showed that infinitely many pairs ($T_{\xi},T_{t}$) existed
such that $T_{\xi}=3T_{t}$. Cunningham \cite{key-3} tested (\ref{eq:1})
for integer values of $k$ (with this paper's notations) by a table
of triangular numbers \cite{key-4,key-6} and found suitable pairs
($\xi,t$) \cite{key-5}. However, solutions can be found also by
posing $x=2\xi+1$ and $y=2t+1$ in (\ref{eq:1}), yielding
\begin{equation}
x^{2}-ky^{2}=\left(1-k\right)\label{eq:2}
\end{equation}
and every odd solutions in ($x,y$) provide the solution pairs ($\xi,t$)
\cite{key-3}. If $k$'s are non-square integers, (\ref{eq:2}) is
a Pell equation that can be solved classically and infinitely many
pairs can be found for each non-square value of $k$ \cite{key-14},
as was rediscovered in \cite{key-7}. Chahal and Top \cite{key-10}
solved the system of two equations $T_{r}=mT_{s}$ and $T_{r}=nT_{t}$
using arithmetic of elliptic curves and theory of elliptic surfaces
to find the infinitely many rational solutions ($r,s,t$).

On the other hand, Gerardin \cite{key-8,key-9} gave a series for
the Diophantine equation $T_{x}\times T_{y}=T_{z}^{2}$ with the recurrence
law $z_{n+1}=6z_{n}-z_{n-1}+2$. Deyi and Tianxin \cite{key-15} proved
that there are infinitely many triangular numbers which can be decomposed
as the product of two triangular numbers, each greater than 1. Pandeyend
\cite{key-18} showed that the product of any two consecutive triangular
numbers can never be a perfect square, but is half another triangular
number; and the product of any two alternate triangular numbers is
the double of another triangular number.

In this paper, we want to find general forms of recursive equations
and closed forms yielding solutions to (\ref{eq:1}). As anticipated,
solutions will be different for $k$ being non-squared integers and
for $k$ squared integers. We are only interested in solutions for
$k>1$ as the cases $k=0$ and $k=1$ are obvious, yielding respectively
$\xi=0$ $\forall t\in\mathbb{Z}^{+}$ and $\xi=t$ $\forall t\in\mathbb{Z}^{+}$.
Although triangular numbers $T_{t}$ are defined for $t\in\mathbb{Z}^{+}$,
we will occasionally consider cases where $t<0$ for which triangular
numbers are defined similarly, yielding $T_{-t}=T_{t-1}$.

\section{Solutions for non-square $k$}

\subsection{Sequences and rank}

Sequences of solutions are known for $k=2,3,5,6,7,8$ and can be found
in the Online Encyclopedia of Integer Sequences (OEIS) \cite{key-19},
as shown in Table 1.

\begin{table}
\caption{OEIS references of sequences of integer solutions of (\ref{eq:1})
for $k=2,3,5,6,7,8$}

\centering{}%
\begin{tabular}{|c||c|c|c|c|}
\hline 
$k$ & $t$ & $\xi$ & $T_{t}$ & $T_{\xi}$\tabularnewline
\hline 
\hline 
2 & A053141 & A001652 & A075528 & A029549\tabularnewline
\hline 
3 & A061278 & A001571 & A076139 & A076140\tabularnewline
\hline 
5 & A077259 & A077262 & A077260 & A077261\tabularnewline
\hline 
6 & A077288 & A077291 & A077289 & A077290\tabularnewline
\hline 
7 & A077398 & A077401 & A077399 & A077400\tabularnewline
\hline 
8 & A336623 & A336625 & A336624 & A336626\tabularnewline
\hline 
\end{tabular}
\end{table}

Let's note first that $t=0$ yielding $\xi=0$ is always a solution
of (\ref{eq:1}) for $\forall k\in\mathbb{Z}$. Let's consider the
two cases of $k=2$ and $k=5$ yielding the successive solution pairs
as shown in Table 2. We indicate also the ratios $t_{n}/t_{n-1}$
for both cases and $t_{n}/t_{n-2}$ for $k=5$. It is seen that for
$k=2,$the ratio $t_{n}/t_{n-1}$ varies between close values, from
7 down to 5.829, while for $k=5$, the ratio $t_{n}/t_{n-1}$ alternates
between values 3 ... 2.619 and 7.333 ... 6.879, while the ratio $t_{n}/t_{n-2}$
decreases regularly from 22 to 18.017 (corresponding approximately
to the product of the alternating values of the ratio $t_{n}/t_{n-1}$).
We call rank $r$ the integer value such that $t_{n}/t_{n-r}$ is
approximately constant or better, decreases regularly without jumps
(a more precise definition is given further). So, here, the case $k=2$
has rank $r=1$ and the case $k=5$ has rank $r=2$.

\begin{table}

\caption{Solutions of (\ref{eq:1}) for $k=2$ and 5}

\centering{}%
\begin{tabular}{|c||c|c|c||c|c|c|c|}
\hline 
$n$ & \multicolumn{3}{c||}{$k=2$} & \multicolumn{4}{c|}{$k=5$}\tabularnewline
\cline{2-8} \cline{3-8} \cline{4-8} \cline{5-8} \cline{6-8} \cline{7-8} \cline{8-8} 
 & $t_{n}$ & $\xi_{n}$ & $t_{n}/t_{n-1}$ & $t_{n}$ & $\xi_{n}$ & $t_{n}/t_{n-1}$ & $t_{n}/t_{n-2}$\tabularnewline
\hline 
\hline 
0 & 0 & 0 &  & 0 & 0 &  & \tabularnewline
\hline 
1 & 2 & 3 & -- & 2 & 5 & -- & --\tabularnewline
\hline 
2 & 14 & 20 & 7 & 6 & 14 & 3 & --\tabularnewline
\hline 
3 & 84 & 119 & 6 & 44 & 99 & 7.33333333 & 22\tabularnewline
\hline 
4 & 492 & 696 & 5.85714286 & 116 & 260 & 2.63636364 & 19.33333333\tabularnewline
\hline 
5 & 2870 & 4059 & 5.83333333 & 798 & 1785 & 6.87931034 & 18.13636364\tabularnewline
\hline 
6 & 16730 & 23660 & 5.82926829 & 2090 & 4674 & 2.61904762 & 18.01724138\tabularnewline
\hline 
\end{tabular}
\end{table}
The rank $r$ is such that the ratio $t_{n}/t_{n-r}$ for $n=2r$,
$t_{2r}/t_{r}$, if corrected by the ratio $t_{r-1}/t_{r}$, is equal
to a constant $2\kappa+3$

\begin{equation}
\frac{t_{2r}-t_{r-1}}{t_{r}}=2\kappa+3\label{eq:3-0}
\end{equation}
where
\begin{equation}
\kappa=t_{r}+t_{r-1}\label{eq:3.1}
\end{equation}
is a constant different for each case of $k.$ For example, for $k=2$
and $r=1$, $\kappa=t_{1}+t_{0}=2$ and for $k=5$ and $r=2$, $\kappa=t_{2}+t_{1}=6+2=8$,
yielding $2\kappa+3=7$ and $19$. This relation (\ref{eq:3-0}) is
important as it links the ratio of successive values of $t_{2r}/t_{r}$
of same rank and the ratio of successive sequential values $t_{r-1}/t_{r}$
to a constant $\kappa$ of the problem.

Another definition of the rank $r$ is as follows. The rank $r$ is
the integer value in the indices of $t_{n}$ and $\xi_{n}$ solutions
of (\ref{eq:1}) such as
\begin{equation}
\kappa=t_{r}+t_{r-1}=\xi_{r}-\xi_{r-1}-1\label{eq:3.2}
\end{equation}
For example, for $k=2$ and $r=1$, $\kappa=\xi_{1}-\xi_{0}-1=3-0-1=2$
and for $k=5$ and $r=2$, $\kappa=\xi_{2}-\xi_{1}-1=14-5-1=8$. This
relation is also important as it relates the two series of basic solutions
$t_{r}$ and $t_{r-1}$on one hand and $\xi_{r}$ and $\xi_{r-1}$
on the other hand.

\subsection{Recurrent relations}

We want to find recurrent relations of the form $X_{n}=C_{1}X_{n-r}+C_{2}X_{n-2r}+C_{3}$
for $t$, $\xi$, $T_{t}$ and $T_{\xi}$ which are solutions of (\ref{eq:1}).
The values of the three coefficients $C_{1},C_{2}$ and $C_{3}$ are
easy to calculate in each case as shown in the following theorem.
\begin{thm}
For $\forall n,k\in\mathbb{Z}^{+}$, $k$ non-square, $\exists r,\kappa,\gamma\in\mathbb{Z}^{+}$
and $\exists t_{n},\xi_{n}\in\mathbb{Z}$ such that there is an infinite
number of solutions of (\ref{eq:1}) of the form

\begin{align}
t_{n} & =2\left(\kappa+1\right)t_{n-r}-t_{n-2r}+\kappa\label{eq:3.3}\\
\xi_{n} & =2\left(\kappa+1\right)\xi_{n-r}-\xi_{n-2r}+\kappa\label{eq:3.3-1}\\
T_{t_{n}} & =\left(4\left(\kappa+1\right)^{2}-2\right)T_{t_{n-r}}-T_{t_{n-2r}}+\left(T_{\kappa}-\gamma\right)\label{eq:3.3-2}\\
T_{\xi_{n}} & =\left(4\left(\kappa+1\right)^{2}-2\right)T_{\xi_{n-r}}-T_{\xi_{n-2r}}+k\left(T_{\kappa}-\gamma\right)\label{eq:3.3-3}
\end{align}
with
\begin{eqnarray}
\kappa & = & t_{r-1}+t_{r}\label{eq:3.4}\\
\gamma & = & t_{r-1}t_{r}\label{eq:3.5}
\end{eqnarray}
\end{thm}

\begin{proof}
Let $n,k,r,\kappa,\gamma,K\in\mathbb{Z}^{+}$, $k$ non-square, and
$t_{n},\xi_{n},C_{1},C_{2},C_{3},C_{4}\in\mathbb{Z}$. Recall first
that $t_{-n}=t_{n-1}$ as $T_{t_{-n}}=T_{t_{n-1}}$, and $t_{0}=0$
is a solution of (\ref{eq:1}), yielding $T_{t_{0}}=0$, $T_{\xi_{0}}=0$,
and $\xi=0$ or $-1$ for $\forall k\in\mathbb{Z}^{+}$, $k$ non-square. 

1) We prove first (\ref{eq:3.3}). Let us write generally $t_{n}=C_{1}t_{n-r}+C_{2}t_{n-2r}+C_{3}$
and assume that $t_{n}$ are solutions of (\ref{eq:1}). Then, for
$n=r-1$, one has $t_{r-1}=C_{1}t_{-1}+C_{2}t_{-r-1}+C_{3}=C_{1}t_{0}+C_{2}t_{r}+C_{3}=C_{2}t_{r}+C_{3}$,
giving
\begin{equation}
C_{3}=t_{r-1}-C_{2}t_{r}\label{eq:4}
\end{equation}
and for $n=r,$$t_{r}=C_{1}t_{0}+C_{2}t_{-r}+C_{3}=C_{2}t_{r-1}+C_{3}$,
yielding
\begin{equation}
C_{3}=t_{r}-C_{2}t_{r-1}\label{eq:5}
\end{equation}
Relations (\ref{eq:4}) and (\ref{eq:5}) are obviously equal, yielding
$C_{2}=-1$ and $C_{3}=t_{r-1}+t_{r}=\kappa$.

For $n=2r,$$t_{2r}=C_{1}t_{r}+C_{2}t_{0}+C_{3}=C_{1}t_{r}+C_{3}=\left(C_{1}+1\right)t_{r}+t_{r-1}$,
yielding $C_{1}=\frac{t_{2r}-t_{r-1}}{t_{r}}-1$. Then, with (\ref{eq:3-0}),
$C_{1}=2\left(\kappa+1\right)$, yielding (\ref{eq:3.3}).

2) Let us write generally $\xi_{n}=C_{1}\xi_{n-r}+C_{2}\xi_{n-2r}+C_{3}$
and assume that $\xi_{n}$ are solutions of (\ref{eq:1}). One has
also that $\xi_{-n}=-\left(\xi_{n-1}+1\right)$ as $T_{\xi_{-n}}=T_{-\left(\xi_{n-1}+1\right)}$
yielding $\xi_{-1}=-1$. For $n=r-1$, one has $\xi_{r-1}=C_{1}\xi_{-1}+C_{2}\xi_{-r-1}+C_{3}$,
giving
\begin{equation}
\xi_{r-1}=-C_{1}-C_{2}\left(\xi_{r}+1\right)+C_{3}\label{eq:9}
\end{equation}
For $n=r$, one has $\xi_{r}=C_{1}\xi_{0}+C_{2}\xi_{-r}+C_{3}=-C_{2}\left(\xi_{r-1}+1\right)+C_{3}$,
yielding
\begin{equation}
\xi_{r}=C_{3}-C_{2}\left(\xi_{r-1}+1\right)\label{eq:10}
\end{equation}
 For $n=2r$, one has $\xi_{2r}=C_{1}\xi_{r}+C_{2}\xi_{0}+C_{3}$,
or
\begin{equation}
\xi_{2r}=C_{1}\xi_{r}+C_{3}\label{eq:11}
\end{equation}
For $n=-1$, one has $\xi_{-1}=C_{1}\xi_{-r-1}+C_{2}\xi_{-2r-1}+C_{3}$
or $-1=-C_{1}\left(\xi_{r}+1\right)-C_{2}\left(\xi_{2r}+1\right)+C_{3}$,
yielding
\begin{equation}
\xi_{2r}=\frac{1-C_{1}\left(\xi_{r}+1\right)-C_{2}+C_{3}}{C_{2}}\label{eq:12}
\end{equation}
Equating (\ref{eq:11}) and (\ref{eq:12}) yields
\begin{equation}
C_{1}=\frac{\left(1-C_{2}\right)\left(C_{3}+1\right)}{\left(C_{2}+1\right)\xi_{r}+1}\label{eq:13}
\end{equation}
Replacing $C_{3}$ from (\ref{eq:10})and $C_{1}$ from (\ref{eq:13})
in (\ref{eq:9}) yield
\begin{equation}
\left(C_{2}-1\right)\left(\xi_{r-1}-\xi_{r}+\frac{C_{2}\left(\xi_{r-1}+1\right)+\xi_{r}+1}{\left(C_{2}+1\right)\xi_{r}+1}\right)=0
\end{equation}
which yields either $C_{2}=1$ or $C_{2}=-1$. The solution with $C_{2}=1$
yields $C_{1}=0$ which is obviously not feasible, and should be discarded.
The other solution $C_{2}=-1$ yields, with (\ref{eq:3.2}), $C_{3}=\xi_{r}-\xi_{r-1}-1=\kappa$
and $C_{1}=2\left(\xi_{r}-\xi_{r-1}\right)=2\left(\kappa+1\right)$,
yielding (\ref{eq:3.3-1}).
\end{proof}
Before proving (\ref{eq:3.3-2}) and (\ref{eq:3.3-3}), we demonstrate
a Lemma that we will need in the proof of (\ref{eq:3.3-2}).
\begin{lem}
For $\forall n,k\in\mathbb{Z}^{+}$, $k$ non-square, $\exists r,\kappa,\gamma\in\mathbb{Z}^{+}$
and $\exists t_{n}\in\mathbb{Z}$ such that $t_{n}$ is solution of
(\ref{eq:1}) and such that the relation
\begin{equation}
\kappa\left(t_{n}+2t_{n}t_{n-r}+t_{n-r}\right)-\left(t_{n}-t_{n-r}\right)^{2}=\gamma\label{eq:5-1}
\end{equation}
is a constant $\gamma=t_{r}t_{r-1}$.
\end{lem}

\begin{proof}
We prove this by induction, showing first that if (\ref{eq:5-1})
holds for $t_{n}$, then (\ref{eq:5-1}) holds also for $t_{n+r}$.
From (\ref{eq:3.3}), one has 
\begin{equation}
t_{n+r}=2\left(\kappa+1\right)t_{n}-t_{n-r}+\kappa\label{eq:5-2}
\end{equation}
Replacing $t_{n}$ and $t_{n-r}$ in (\ref{eq:5-1}) by $t_{n+r}$
and $t_{n}$ and then $t_{n+r}$ by (\ref{eq:5-2}) yield directly
\begin{eqnarray}
 &  & \kappa\left(t_{n+r}+2t_{n+r}t_{n}+t_{n}\right)-\left(t_{n+r}-t_{n}\right)^{2}=\nonumber \\
 &  & \kappa\left(t_{n}+2t_{n}t_{n-r}+t_{n-r}\right)-\left(t_{n}-t_{n-r}\right)^{2}=\gamma\label{eq:5-3}
\end{eqnarray}
As $r\geq1$ and as $n$ can take any positive integer value, we have
showed that (\ref{eq:5-1}) for $t_{n}$ implies directly (\ref{eq:5-3})
for $t_{n+r}$.

For $n=0$, (\ref{eq:5-1}) yields, with (\ref{eq:3.1}), $\gamma=\kappa\left(t_{0}+2t_{0}t_{-r}+t_{-r}\right)-\left(t_{0}-t_{-r}\right)^{2}=\kappa t_{r-1}-t_{r-1}^{2}=\left(t_{r}+t_{r-1}\right)t_{r-1}-t_{r-1}^{2}=t_{r}t_{r-1}$
\end{proof}
We proceed now with the rest of the proof of Theorem 1.
\begin{proof}
3) Calculating $T_{t_{n}}=\frac{t_{n}^{2}+t_{n}}{2}$ with (\ref{eq:3.3})
yields successively
\begin{eqnarray}
T_{t_{n}} & = & \frac{1}{2}\left[\left(2\left(\kappa+1\right)t_{n-r}-t_{n-2r}+\kappa\right)^{2}+\left(2\left(\kappa+1\right)t_{n-r}-t_{n-2r}+\kappa\right)\right]\nonumber \\
 & = & 4\left(\kappa+1\right)^{2}T_{t_{n-r}}+T_{t_{n-2r}}+T_{\kappa}-\left(\kappa+1\right)\left(t_{n-r}+2t_{n-r}t_{n-2r}+t_{n-2r}\right)\nonumber \\
 & = & \left(4\left(\kappa+1\right)^{2}-2\right)T_{t_{n-r}}-T_{t_{n-2r}}+C_{4}\label{eq:6}
\end{eqnarray}
which is (\ref{eq:3.3-2}) and where $C_{4}$ is a constant by Lemma
2,
\begin{eqnarray}
C_{4} & = & T_{\kappa}+2\left(T_{t_{n-r}}+T_{t_{n-2r}}\right)-\left(\kappa+1\right)\left(t_{n-r}+2t_{n-r}t_{n-2r}+t_{n-2r}\right)\nonumber \\
 & = & T_{\kappa}-\kappa\left(t_{n-r}+2t_{n-r}t_{n-2r}+t_{n-2r}\right)+\left(t_{n-r}-t_{n-2r}\right)^{2}=T_{\kappa}-\gamma\label{eq:7}
\end{eqnarray}
4) Replacing $T_{t_{n}}$in (\ref{eq:1}) by (\ref{eq:3.3-2}) yields
\begin{eqnarray}
T_{\xi_{n}} & = & \left(4\left(\kappa+1\right)^{2}-2\right)kT_{t_{n-r}}-kT_{t_{n-2r}}+k\left(T_{\kappa}-\gamma\right)\nonumber \\
 & = & \left(4\left(\kappa+1\right)^{2}-2\right)T_{\xi_{n-r}}-T_{\xi_{n-2r}}+k\left(T_{\kappa}-\gamma\right)\label{eq:8}
\end{eqnarray}
which is (\ref{eq:3.3-3}) and where (\ref{eq:1}) was applied to
$T_{\xi_{n-r}}$ and $T_{\xi_{n-2r}}$.

5) As $n$ is unbound, the number of solutions given by (\ref{eq:3.3})
to (\ref{eq:3.3-3}) is infinite.
\end{proof}
Note that relations (\ref{eq:3.3}) to (\ref{eq:3.3-3}) are independent
from the value of $k$. Only the rank $r$ depends in an unknown manner
on $k$. With these relations (\ref{eq:3.3}) to (\ref{eq:3.3-3}),
one can easily calculate all solutions of (\ref{eq:1}) provided that
the rank $r$ is known and the $n$ first values of $t_{n}$ are known
for $1\leq n\leq r$. 

Table 3 shows the first six sets of recurrent equations for $t_{n},\xi_{n},T_{t_{n}},T_{\xi_{n}}$
for non-square $k$, $2\leq k\leq8$, with the value of $r$ and the
first $t_{2r}$ values (from references in Table 1).

\begin{table}

\begin{raggedright}
\caption{Recurrent equations for $t_{n},\xi_{n},T_{t_{n}},T_{\xi_{n}}$ for
non-square $k$, $2\protect\leq k\protect\leq8$}
\par\end{raggedright}
\begin{raggedright}
\uline{\mbox{$k=2,r=1$}}
\par\end{raggedright}
\begin{raggedright}
$t_{n}=6t_{n-1}-t_{n-2}+2$, $t_{n}=0,2$
\par\end{raggedright}
\begin{raggedright}
$\xi_{n}=6\xi_{n-1}-\xi_{n-2}+2$, $\xi_{n}=0,3$
\par\end{raggedright}
\begin{raggedright}
$Tt_{n}=34Tt_{n-1}-Tt_{n-2}+3$,$Tt_{n}=0,3$
\par\end{raggedright}
\begin{raggedright}
$T\xi_{n}=34T\xi_{n-1}-T\xi_{n-2}+6$, $T\xi_{n}=0,6$
\par\end{raggedright}
\begin{raggedright}
\uline{\mbox{$k=3,r=1$}}
\par\end{raggedright}
\begin{raggedright}
$t_{n}=4t_{n-1}-t_{n-2}+1$, $t_{n}=0,1$
\par\end{raggedright}
\begin{raggedright}
$\xi_{n}=4\xi_{n-1}-\xi_{n-2}+1$, $\xi_{n}=0,2$
\par\end{raggedright}
\begin{raggedright}
$Tt_{n}=14Tt_{n-1}-Tt_{n-2}+1$, $Tt_{n}=0,1$
\par\end{raggedright}
\begin{raggedright}
$T\xi_{n}=14T\xi_{n-1}-T\xi_{n-2}+3$, $T\xi_{n}=0,3$
\par\end{raggedright}
\begin{raggedright}
\uline{\mbox{$k=5,r=2$}}
\par\end{raggedright}
\begin{raggedright}
$t_{n}=18t_{n-2}-t_{n-4}+8$, $t_{n}=0,2,6,44$
\par\end{raggedright}
\begin{raggedright}
$\xi_{n}=18\xi_{n-2}-\xi_{n-4}+8$, $\xi_{n}=0,5,14,99$
\par\end{raggedright}
\begin{raggedright}
$Tt_{n}=322Tt_{n-2}-Tt_{n-4}+24$, $Tt_{n}=0,3,21,990$
\par\end{raggedright}
\begin{raggedright}
$T\xi_{n}=322T\xi_{n-2}-T\xi_{n-4}+120$, $T\xi_{n}=0,15,105,4950$
\par\end{raggedright}
\begin{raggedright}
\uline{\mbox{$k=6,r=2$}}
\par\end{raggedright}
\begin{raggedright}
$t_{n}=10t_{n-2}{}_{-}t_{n-4}{}_{+}4,t_{n}=0,1,3,14$
\par\end{raggedright}
\begin{raggedright}
$\xi_{n}=10\xi_{n-2}-\xi_{n-4}+4,\xi_{n}=0,3,8,35$
\par\end{raggedright}
\begin{raggedright}
$Tt_{n}=98Tt_{n-2}-Tt_{n-4}+7,Tt_{n}=0,1,6,105$
\par\end{raggedright}
\begin{raggedright}
$T\xi_{n}=98T\xi_{n-2}-T\xi_{n-4}+42,T\xi_{n}=0,6,36,630$
\par\end{raggedright}
\begin{raggedright}
\uline{\mbox{$k=7,r=2$}}
\par\end{raggedright}
\begin{raggedright}
$t_{n}=16t_{n-2}-t_{n-4}+7,t_{n}=0,2,5,39$
\par\end{raggedright}
\begin{raggedright}
$\xi_{n}=16\xi_{n-2}-\xi_{n-4}+7,\xi_{n}=0,6,14,104$
\par\end{raggedright}
\begin{raggedright}
$Tt_{n}=254Tt_{n-2}-Tt_{n-4}+18,Tt_{n}=0,3,15,780$
\par\end{raggedright}
\begin{raggedright}
$T\xi_{n}=254T\xi_{n-2}-T\xi_{n-4}+126,T\xi_{n}=0,21,105,5460$
\par\end{raggedright}
\begin{raggedright}
\uline{\mbox{$k=8,r=2$}}
\par\end{raggedright}
\begin{raggedright}
$t_{n}=34t_{n-2}-t_{n-4}+16,t_{n}=0,5,11,186$
\par\end{raggedright}
\begin{raggedright}
$\xi_{n}=34\xi_{n-2}-\xi_{n-4}+16,\xi_{n}=0,15,32,527$
\par\end{raggedright}
\begin{raggedright}
$Tt_{n}=1154Tt_{n-2}-Tt_{n-4}+81,Tt_{n}=0,15,66,17391$
\par\end{raggedright}
\raggedright{}$T\xi_{n}=1154T\xi_{n-2}-T\xi_{n-4}+648,T\xi_{n}=0,120,528,139128$
\end{table}

\subsection{Other relations}

There are several other interesting relations between $t_{n}$, $\xi_{n}$,
$T_{t_{n}}$ and $T_{\xi_{n}}$ as follows, which are all valid $\forall n,r,k\in\mathbb{Z}^{+}$,
\begin{align}
t_{n-r}t_{n+r} & =\left(t_{n}-t_{r-1}\right)\left(t_{n}-t_{r}\right)\label{eq:14}\\
\xi_{n-r}\xi_{n+r} & =\left(\xi_{n}+1-\xi_{r-1}\right)\left(\xi_{n}-\xi_{r}\right)\label{eq:15}\\
T_{t_{n-r}}T_{t_{n+r}} & =\left(T_{t_{n}}-T_{t_{r-1}}\right)\left(T_{t_{n}}-T_{t_{r}}\right)\label{eq:16}\\
T_{\xi_{n-r}}T_{\xi_{n+r}} & =\left(T_{\xi_{n}}-T_{\xi_{r-1}}\right)\left(T_{\xi_{n}}-T_{\xi_{r}}\right)\label{eq:17}
\end{align}
and
\begin{equation}
\left(2t_{r-1}+1\right)\left(t_{r}+t_{r-1}\right)=t_{2r-1}-2t_{r-1}\label{eq:18}
\end{equation}
\begin{equation}
\frac{t_{r}}{t_{r-1}}=\frac{t_{2r}-t_{r-1}}{t_{2r-1}-t_{r}}\label{eq:19}
\end{equation}
There is also a relation in $\xi_{n}$ similar to (\ref{eq:5-1})
for $t_{n}$
\begin{equation}
\left(\xi_{n}-\xi_{n-r}\right)^{2}-\kappa\left(\xi_{n}+2\xi_{n}\xi_{n-r}+\xi_{n-r}\right)=\xi_{r}\left(\xi_{r-1}+1\right)\label{eq:20}
\end{equation}
and that can be demonstrated by induction like in Lemma 2. Very interestingly,
the relation (\ref{eq:3.2}) can be generalized as
\begin{equation}
t_{\mu r}+t_{\mu r-1}=\xi_{\mu r}-\xi_{\mu r-1}-1\label{eq:21}
\end{equation}
where the rank $r$ index of $t$ and $\xi$ is replaced by a multiple
$\mu r$ with $\forall\mu\in\mathbb{Z}^{+}$. 

For all the $2s$ values of $k$ between two successive squares, $s^{2}<k<s^{\prime2}=\left(s+1\right)^{2}$,
the ratios $t_{r}/t_{r-1}$ are such that $\left(s+1\right)/\left(s-1\right)\leq t_{r}/t_{r-1}<\left(s^{\prime}+1\right)/\left(s^{\prime}-1\right)=\left(s+2\right)/s$,
i.e. the $2s$ ratios $t_{r}/t_{r-1}$ always vary between these limits.
This means also that for $s>>1$, $\left(s+1\right)/\left(s-1\right)=\left(1+1/s\right)/\left(1-1/s\right)\approx\left(1+1/s\right)^{2}\approx\left(1+2/s\right)$,
and $t_{r}/t_{r-1}\approx\left(1+2/s\right)$

And finally and obviously, the ratio $\frac{\xi_{n}}{t_{n}}$ is such
that $\lim_{n\rightarrow\infty}\left(\frac{\xi_{n}}{t_{n}}\right)=\sqrt{k}$
, as $k=\frac{T_{\xi_{n}}}{T_{t_{n}}}=\frac{\xi_{n}\left(\xi_{n}+1\right)}{t_{n}\left(t_{n}+1\right)}\approx\frac{\xi_{n}^{2}}{t_{n}^{2}}$
for $t_{n}>>1$ and $\xi_{n}>>1$.

\subsection{Values of $\kappa$ in function of $k$}

Recursive relations still depends on the parameters $\kappa$, $\gamma$
and rank $r$. Table 4 gives the values of $r$, $t_{r-1}$, $t_{r}$,
$\kappa$ for non-square values of $k$ between 2 and 102, found numerically.

\begin{table}
\begin{centering}
\caption{Values of $k$, $r,$$t_{r-1}$, $t_{r}$ and $\kappa$}
\par\end{centering}
\centering{}%
\begin{tabular}{|c|c|c|c|c||c|c|c|c|c|}
\hline 
$k$ & $r$ & $t_{r-1}$ & $t_{r}$ & $\kappa$  & $k$ & $r$ & $t_{r-1}$ & $t_{r}$ & $\kappa$ \tabularnewline
\hline 
\hline 
2 & 1 & 0 & 2 & 2 & 54 & 2 & 209 & 275 & 484\tabularnewline
\hline 
3 & 1 & 0 & 1 & 1 & 55 & 4 & 38 & 50 & 88\tabularnewline
\hline 
5 & 2 & 2 & 6 & 8 & 56 & 2 & 6 & 8 & 14\tabularnewline
\hline 
6 & 2 & 1 & 3 & 4 & 57 & 4 & 65 & 85 & 150\tabularnewline
\hline 
7 & 2 & 2 & 5 & 7 & 58 & 4 & 8514 & 11088 & 19602\tabularnewline
\hline 
8 & 2 & 5 & 11 & 16 & 59 & 2 & 230 & 299 & 529\tabularnewline
\hline 
10 & 3 & 6 & 12 & 18 & 60 & 2 & 13 & 17 & 30\tabularnewline
\hline 
11 & 2 & 3 & 6 & 9 & 61 & 8 & 770082534 & 996236514 & 1766319048\tabularnewline
\hline 
12 & 2 & 2 & 4 & 6 & 62 & 2 & 27 & 35 & 62\tabularnewline
\hline 
13 & 4 & 234 & 414 & 648 & 63 & 2 & 3 & 4 & 7\tabularnewline
\hline 
14 & 2 & 5 & 9 & 14 & 65 & 2 & 56 & 72 & 128\tabularnewline
\hline 
15 & 2 & 1 & 2 & 3 & 66 & 4 & 28 & 36 & 64\tabularnewline
\hline 
17 & 2 & 12 & 20 & 32 & 67 & 4 & 21437 & 27404 & 48841\tabularnewline
\hline 
18 & 2 & 6 & 10 & 16 & 68 & 2 & 14 & 18 & 32\tabularnewline
\hline 
19 & 3 & 65 & 104 & 169 & 69 & 4 & 3419 & 4355 & 7774\tabularnewline
\hline 
20 & 2 & 3 & 5 & 8 & 70 & 4 & 110 & 140 & 250\tabularnewline
\hline 
21 & 4 & 21 & 33 & 54 & 71 & 4 & 1533 & 1946 & 3479\tabularnewline
\hline 
22 & 4 & 77 & 119 & 196 & 72 & 2 & 7 & 9 & 16\tabularnewline
\hline 
23 & 2 & 9 & 14 & 23 & 73 & 6 & 1007124 & 1274124 & 2281248\tabularnewline
\hline 
24 & 2 & 19 & 29 & 48 & 74 & 2 & 1634 & 2064 & 3698\tabularnewline
\hline 
26 & 3 & 20 & 30 & 50 & 75 & 2 & 11 & 14 & 25\tabularnewline
\hline 
27 & 2 & 10 & 15 & 25 & 76 & 6 & 25584 & 32214 & 57798\tabularnewline
\hline 
28 & 4 & 51 & 75 & 126 & 77 & 4 & 155 & 195 & 350\tabularnewline
\hline 
29 & 4 & 3990 & 5810 & 9800 & 78 & 4 & 23 & 29 & 52\tabularnewline
\hline 
30 & 2 & 4 & 6 & 10 & 79 & 2 & 35 & 44 & 79\tabularnewline
\hline 
31 & 4 & 623 & 896 & 1519 & 80 & 2 & 71 & 89 & 160\tabularnewline
\hline 
32 & 2 & 237 & 339 & 576 & 82 & 3 & 72 & 90 & 162\tabularnewline
\hline 
33 & 2 & 9 & 13 & 22 & 83 & 2 & 36 & 45 & 81\tabularnewline
\hline 
34 & 2 & 14 & 20 & 34 & 84 & 2 & 24 & 30 & 54\tabularnewline
\hline 
35 & 2 & 2 & 3 & 5 & 85 & 8 & 127386 & 158382 & 285768\tabularnewline
\hline 
37 & 2 & 30 & 42 & 72 & 86 & 4 & 4641 & 5763 & 10404\tabularnewline
\hline 
38 & 2 & 15 & 21 & 36 & 87 & 2 & 12 & 15 & 27\tabularnewline
\hline 
39 & 2 & 10 & 14 & 24 & 88 & 4 & 34671 & 42945 & 77616\tabularnewline
\hline 
40 & 4 & 303 & 417 & 720 & 89 & 4 & 223500 & 276500 & 500000\tabularnewline
\hline 
41 & 4 & 864 & 1184 & 2048 & 90 & 2 & 8 & 10 & 18\tabularnewline
\hline 
42 & 2 & 5 & 7 & 12 & 91 & 6 & 704 & 869 & 1573\tabularnewline
\hline 
43 & 4 & 1475 & 2006 & 3481 & 92 & 4 & 515 & 635 & 1150\tabularnewline
\hline 
44 & 2 & 84 & 114 & 198 & 93 & 4 & 5445 & 6705 & 12150\tabularnewline
\hline 
45 & 4 & 68 & 92 & 160 & 94 & 4 & 961115 & 1182179 & 2143294\tabularnewline
\hline 
46 & 6 & 10373 & 13961 & 24334 & 95 & 2 & 17 & 21 & 38\tabularnewline
\hline 
47 & 2 & 20 & 27 & 47 & 96 & 4 & 2155 & 2645 & 4800\tabularnewline
\hline 
48 & 2 & 41 & 55 & 96 & 97 & 4 & 28216140 & 34593492 & 62809632\tabularnewline
\hline 
50 & 3 & 42 & 56 & 98 & 98 & 2 & 44 & 54 & 98\tabularnewline
\hline 
51 & 3 & 21 & 28 & 49 & 99 & 2 & 4 & 5 & 9\tabularnewline
\hline 
52 & 4 & 279 & 369 & 648 & 101 & 2 & 90 & 110 & 200\tabularnewline
\hline 
53 & 4 & 28574 & 37674 & 66248 & 102 & 2 & 45 & 55 & 100\tabularnewline
\hline 
\end{tabular}
\end{table}

\section{Solutions for square $k$}

\subsection{Recurrent relations}

Contrary to the cases of $k$ not being a square, for $k$ being a
square, there are no infinitely many solutions. Depending on the values
of $k,$ there are no, or only one, solution.

Let us write $k=\lambda^{2}$, with $\lambda\in\mathbb{Z}$. Relation
(\ref{eq:1}) reads then
\begin{equation}
T_{\xi}=\lambda^{2}T_{t}\label{eq:1-1}
\end{equation}

As said previously, $\lambda=0$ yields $\xi=0$, $\forall t$, while
$\lambda=1$ yields $\xi=t$, $\forall t$, and are therefore always
solution of (\ref{eq:1-1}). Other solutions are given in the following
Lemma where we search for the values of $\lambda$ and $\xi$ yielding
solutions to (\ref{eq:1-1}) for single values of $t$.
\begin{lem}
For $\forall t,n\in\mathbb{Z^{+}}$, $\exists\lambda,\xi\in\mathbb{Z}$
such that $\lambda_{n},\xi_{n}$ are solutions of (\ref{eq:1-1})
and $\lambda_{n}$ and $\xi_{n}$ are polynomials of increasing degrees
in $t$, $\lambda_{n}=P(t)_{\left(n-1\right)}$ of $\left(n-1\right)$degree
and $\xi_{n}=Q(t)_{\left(n\right)}$ of $n$ degree, whose smallest
are the polynomials for $n=2$
\begin{align}
\lambda_{2} & =4t+2\label{eq:3.2-1}\\
\xi_{2} & =4t\left(t+1\right)=8T_{t}\label{eq:3.2.2}
\end{align}
\end{lem}

\begin{proof}
Let $t\in\mathbb{Z^{+}}$, and $\lambda,\xi,A,B,C,D,E\in\mathbb{Z}$.
For (\ref{eq:1-1}) to hold, $\lambda^{2}t\left(t+1\right)/2$ must
be a triangular number or $4\lambda^{2}t\left(t+1\right)+1$ must
be a square. Let $\lambda$ be a first degree polynomial $\lambda=At+B$
where $A$ and $B$ are unknown integer constants. Then $4\left(At+B\right)^{2}\left(t^{2}+t\right)+1$
must be a squared second degree polynomial, namely $4\left(At+B\right)^{2}\left(t^{2}+t\right)+1=\left(Ct^{2}+Dt+E\right)^{2}$
where $C$, $D$ and $E$ are unknown integer constants. Calculating
and equating coefficients of same power of $t$ yield five conditions
to determine the five constants. One has $C=\pm2A$, $CD=2A\left(A+2B\right)$,
$2CE+D^{2}=4B\left(2A+B\right)$, $DE=2B^{2}$ and $E=\pm1$. Solving
with simple algebra yields eventually $A=4$, $B=2$, $C=8$, $D=8$
and $E=1$, yielding $4\left(4t+2\right)^{2}\left(t^{2}+t\right)+1=\left(8t^{2}+8t+1\right)^{2}$
and $\lambda^{2}t\left(t+1\right)/2=\left(4t+2\right)^{2}\left(t^{2}+t\right)/2=\left(4t^{2}+4t+1\right)\left(4t^{2}+4t\right)/2=T_{\left(4t^{2}+4t\right)}=T_{\xi}$
yielding $\xi=\left(4t^{2}+4t\right)=8T_{t}$.
\end{proof}
Relation (\ref{eq:3.2-1}) shows that there is a solution for each
integer value of $t$ to which corresponds a single value of $\lambda$
(or $k$) and of $\xi$. Other polynomials $\lambda_{n}=P(t)_{\left(n-1\right)}$
of higher degrees $\left(n-1\right)$ and $\xi_{n}=Q(t)_{\left(n\right)}$
of higher degrees $n$, can be calculated similarly as shown in Table
5 for the first few polynomials $\lambda_{n}=P(t)_{\left(n-1\right)}$.

\begin{table}

\caption{Polynomials $\lambda_{n}=P(t)_{\left(n-1\right)}$}

\centering{}%
\begin{tabular}{|c|l|}
\hline 
$n$ & $\lambda_{n}$\tabularnewline
\hline 
\hline 
2 & $4t+2$\tabularnewline
\hline 
3 & $\left(4t+2\right)^{2}-1$\tabularnewline
\hline 
4 & $\left[\left(4t+2\right)^{2}-2\right]\left(4t+2\right)$\tabularnewline
\hline 
5 & $\left[\left(4t+2\right)^{2}-3\right]\left(4t+2\right)^{2}+1$\tabularnewline
\hline 
6 & $\left\{ \left[\left(4t+2\right)^{2}-4\right]\left(4t+2\right)^{2}+3\right\} \left(4t+2\right)$\tabularnewline
\hline 
7 & $\left\{ \left[\left(4t+2\right)^{2}-5\right]\left(4t+2\right)^{2}+6\right\} \left(4t+2\right)^{2}-1$\tabularnewline
\hline 
8 & $\left\{ \left[\left(4t+2\right)^{2}-6\right]\left(4t+2\right)^{2}+10\right\} \left[\left(4t+2\right)^{2}-4\right]\left(4t+2\right))$\tabularnewline
\hline 
\end{tabular}
\end{table}

However, this process becomes more and more tedious with increasing
$n$. A simpler approach is to use a recurrent relation as given in
the following theorem.
\begin{thm}
For $\forall t,n\in\mathbb{Z^{+}}$, $\exists\lambda,\xi\in\mathbb{Z}$
such that $\lambda,\xi$ are solutions of (\ref{eq:1-1}) of the form
\begin{align}
\lambda_{n} & =\left(4t+2\right)\lambda_{n-1}-\lambda_{n-2}\label{eq:3.4-1}\\
\xi_{n} & =\left(4t+2\right)\xi_{n-1}-\xi_{n-2}+2t\label{eq:3.5-1}
\end{align}
with $\lambda_{-1}=-1$, $\lambda_{0}=0$, $\lambda_{1}=1$, $\xi_{0}=0$
and $\xi_{1}=\xi_{-1}=t$.
\end{thm}

\begin{proof}
Let $t,\xi\in\mathbb{Z^{+}}$ and $n,\lambda,C_{1},C_{2},C_{3}\in\mathbb{Z}$.

1) We prove first (\ref{eq:3.4-1}). Let us write generally $\lambda_{n}=C_{1}\lambda_{n-1}+C_{2}\lambda_{n-2}+C_{3}$
and assume that $\lambda_{n}$ are solutions of (\ref{eq:1-1})with
$\lambda_{-1}=-1$, $\lambda_{0}=0$ and $\lambda_{1}=1$. By Lemma
3, we know also that $\lambda_{2}=4t+2$ is a solution of (\ref{eq:1-1}),
yielding the additional condition $\lambda_{-2}=-\left(4t+2\right)$.
Writing (\ref{eq:3.4-1}) for $n=0$ to $2$ yields successively
\begin{eqnarray}
\lambda_{0}=C_{1}\lambda_{-1}+C_{2}\lambda_{-2}+C_{3} & \Rightarrow & 0=-C_{1}-C_{2}\left(4t+2\right)+C_{3}\label{eq:3.6}\\
\lambda_{1}=C_{1}\lambda_{0}+C_{2}\lambda_{-1}+C_{3} & \Rightarrow & 1=-C_{2}+C_{3}\label{eq:3.7}\\
\lambda_{2}=C_{1}\lambda_{1}+C_{2}\lambda_{0}+C_{3} & \Rightarrow & 4t+2=C_{1}+C_{3}\label{eq:3.8}
\end{eqnarray}
Solving the three equations (\ref{eq:3.6}) to (\ref{eq:3.8}) for
$C_{1},C_{2},C_{3}$, one finds $C_{1}=4t+2$, $C_{2}=-1$ and $C_{3}=0$,
yielding (\ref{eq:3.4-1}).

2) For (\ref{eq:3.5-1}), let us write generally $\xi_{n}=C_{1}\xi_{n-1}+C_{2}\xi_{n-2}+C_{3}$
and assume that $\xi_{n}$ are solutions of (\ref{eq:1-1}) with $\xi_{0}=0$
(for $\lambda_{0}=0$) and $\xi_{1}=\xi_{-1}=t$ (for $\lambda_{\pm1}=\pm1$).
By Lemma 3, we know also that $\xi_{2}=4t\left(t+1\right)$ is a solution
of (\ref{eq:1-1}), yielding the additional condition $\xi_{-2}=4t\left(t+1\right)$.
Writing (\ref{eq:3.5-1}) for $n=0$ to $2$ yields successively
\begin{eqnarray}
\xi_{0}=C_{1}\xi_{-1}+C_{2}\xi_{-2}+C_{3} & \Rightarrow & 0=C_{1}t+C_{2}4t\left(t+1\right)+C_{3}\label{eq:3.9}\\
\xi_{1}=C_{1}\xi_{0}+C_{2}\xi_{-1}+C_{3} & \Rightarrow & t=C_{2}t+C_{3}\label{eq:3.10}\\
\xi_{2}=C_{1}\xi_{1}+C_{2}\xi_{0}+C_{3} & \Rightarrow & 4t\left(t+1\right)=C_{1}t+C_{3}\label{eq:3.11}
\end{eqnarray}
Solving the three equations (\ref{eq:3.9}) to (\ref{eq:3.11}) for
$C_{1},C_{2},C_{3}$ yields $C_{1}=4t+2$, $C_{2}=-1$ and $C_{3}=2t$,
yielding (\ref{eq:3.5-1}).
\end{proof}

\subsection{Values of $k$ and $\xi$}

Table 6 shows the first $n=5$ values of $k_{n}=\lambda_{n}^{2}$
and $\xi_{n}$ from (\ref{eq:3.4-1}) and (\ref{eq:3.5-1}) for $1\leq t\leq6$.

\begin{table}
\caption{Values of $k_{n}=\lambda_{n}^{2}$ and $\xi_{n}$}

\centering{}%
\begin{tabular}{|c||c|c||c|c||c|c|}
\hline 
 & \multicolumn{2}{c||}{$t=1$} & \multicolumn{2}{c||}{$t=2$} & \multicolumn{2}{c|}{$t=3$}\tabularnewline
\cline{2-7} \cline{3-7} \cline{4-7} \cline{5-7} \cline{6-7} \cline{7-7} 
$n$ & $k_{n}$ & $\xi_{n}$ & $k_{n}$ & $\xi_{n}$ & $k_{n}$ & $\xi_{n}$\tabularnewline
\hline 
\hline 
1 & 1 & 1 & 1 & 2 & 1 & 3\tabularnewline
\hline 
2 & 36 & 8 & 100 & 24 & 196 & 48\tabularnewline
\hline 
3 & 1225 & 49 & 9801 & 242 & 38025 & 675\tabularnewline
\hline 
4 & 41616 & 288 & 960400 & 2400 & 7376656 & 9408\tabularnewline
\hline 
5 & 1413721 & 1681 & 94109401 & 23762 & 1431033241 & 131043\tabularnewline
\hline 
\hline 
 & \multicolumn{2}{c||}{$t=4$} & \multicolumn{2}{c||}{$t=5$} & \multicolumn{2}{c|}{$t=6$}\tabularnewline
\cline{2-7} \cline{3-7} \cline{4-7} \cline{5-7} \cline{6-7} \cline{7-7} 
$n$ & $k_{n}$ & $\xi_{n}$ & $k_{n}$ & $\xi_{n}$ & $k_{n}$ & $\xi_{n}$\tabularnewline
\hline 
1 & 1 & 4 & 1 & 5 & 1 & 6\tabularnewline
\hline 
2 & 324 & 80 & 484 & 120 & 676 & 168\tabularnewline
\hline 
3 & 104329 & 1444 & 233289 & 2645 & 455625 & 4374\tabularnewline
\hline 
4 & 33593616 & 25920 & 112444816 & 58080 & 307090576 & 113568\tabularnewline
\hline 
5 & 10817040025 & 465124 & 54198168025 & 1275125 & 206978592601 & 2948406\tabularnewline
\hline 
\end{tabular}
\end{table}

One notices that the value of $k_{n}$ for $t=1$ are all squares
and triangular numbers (see OEIS, A001110) which is obvious as if
$t=1$, $T_{t}=1$ and for (\ref{eq:1-1}) to hold, $k$ must be also
a triangular number in addition of being a square. Note that the $\xi_{n}$
for $t=1$ are given in OEIS, A001108. 

Table 7 lists the values of $\lambda,k,t$ and $\xi$ for sequentially
increasing values of $\lambda$ and $k<10^{4}$.

\begin{table}

\caption{Values of $\lambda,k,t$ and $\xi$}

\centering{}%
\begin{tabular}{|c|c|c|c||c|c|c|c||c|c|c|c|}
\hline 
$\lambda$ & $k$ & $t$ & $\xi$ & $\lambda$ & $k$ & $t$ & $\xi$ & $\lambda$ & $k$ & $t$ & $\xi$\tabularnewline
\hline 
\hline 
6 & 36 & 1 & 8 & 38 & 1444 & 9 & 360 & 74 & 5476 & 18 & 1368\tabularnewline
\hline 
10 & 100 & 2 & 24 & 42 & 1764 & 10 & 440 & 78 & 6084 & 19 & 1520\tabularnewline
\hline 
14 & 196 & 3 & 48 & 46 & 2116 & 11 & 528 & 82 & 6724 & 20 & 1680\tabularnewline
\hline 
18 & 324 & 4 & 80 & 50 & 2500 & 12 & 624 & 86 & 7396 & 21 & 1848\tabularnewline
\hline 
22 & 484 & 5 & 120 & 54 & 2916 & 13 & 728 & 90 & 8100 & 22 & 2024\tabularnewline
\hline 
26 & 676 & 6 & 168 & 58 & 3364 & 14 & 840 & 94 & 8836 & 23 & 2208\tabularnewline
\hline 
30 & 900 & 7 & 224 & 62 & 3844 & 15 & 960 & 98 & 9604 & 24 & 2400\tabularnewline
\hline 
34 & 1156 & 8 & 288 & 66 & 4356 & 16 & 1088 & 99 & 9801 & 2 & 242\tabularnewline
\hline 
35 & 1225 & 1 & 49 & 70 & 4900 & 17 & 1224 &  &  &  & \tabularnewline
\hline 
\end{tabular}
\end{table}

\section{Conclusions}

We have shown that for any positive non-square integer multiplier
$k$, there is an infinity of multiples of triangular numbers that
are triangular numbers. Recurrent relations were deduced. For squared
integer values of the multiplier $k$, there is either one or no solution,
depending on the value of $k$ as shown in Table 7. In a following
paper, we explore other properties of triangular numbers that are
multiple of other triangular numbers.

\end{document}